\documentclass[12pt]{article}
\usepackage{amsmath,amssymb}

\newcommand{\ncm}{\newcommand}
\ncm{\aut}{auto\-mor\-phi\-sm} \ncm{\Inn}{\mbox{\rm Inn}} 
\ncm{\Ap}{\mbox{$\overline{\rm Inn}$}} \ncm{\Ext}{\mbox{\rm Ext}} 
\ncm{\Ex}{\mbox{\rm Ex}} \ncm{\OExt}{\mbox{\rm OrderExt}} 
\ncm{\AI}{\mbox{\rm AInn}} \ncm{\HI}{\mbox{\rm HInn($A$)}} 
\ncm{\Aut}{\mbox{\rm Aut}} \ncm{\Mal}{\mbox{$M_{\alpha}$}} 
\ncm{\Aff}{\mbox{${\rm Aff}$}} \ncm{\id}{\mbox{\rm id}} 
\ncm{\Ker}{\mbox{\rm Ker}} \ncm{\BE}{\begin{eqnarray*}} 
\ncm{\EE}{\end{eqnarray*}} \ncm{\lra}{\mbox{$\longrightarrow$}} 
\ncm{\Hom}{\mbox{\rm Hom}} \ncm{\calU}{{\cal U}} \ncm{\el}{\ell} 
\ncm{\ad}{\mbox{\rm ad}} \ncm{\Alg}{\mbox{\rm Alg}} 
\ncm{\Conv}{\mbox{\rm Conv}} \ncm{\D}{{\cal D}} 

\ncm{\cstar}{$C^{*}$-algebra} \ncm{\cstars}{$C^{*}$-algebras} 
\ncm{\ra}{\mbox{$\rightarrow$}} \ncm{\la}{\mbox{$\leftarrow$}} 
\ncm{\hra}{\hookrightarrow} \ncm{\da}{\mbox{$\downarrow$}} 
\ncm{\se}{\mbox{$\searrow$}} \ncm{\del}{\mbox{$\delta$}} 
\ncm{\supp}{\mbox{\rm supp}} \ncm{\Ad}{\mbox{\rm Ad}} 
\ncm{\CAR}{\mbox{$M_{2^{\infty}}$}} \ncm{\ep}{\mbox{$\epsilon > 
0$}} \ncm{\ol}{\overline} \ncm{\Mninf}{\mbox{$M_{n^{\infty}}$}} 
\ncm{\MR}{M. R\o{}rdam} \ncm{\Range}{\mbox{\rm Range}} 
\ncm{\vo}{}
\ncm{\ch}{}
\ncm{\CMP}{Comm. Math. Phys.} \ncm{\add}{} 
\ncm{\tilsig}{\tilde{\sigma}} \ncm{\dist}{{\rm 
dist}}\ncm{\eps}{\epsilon}  \ncm{\calL}{{\mathcal{L}}} 
 
\ncm{\calH}{{\mathcal{H}}} 
\ncm{\lan}{{\langle}}\ncm{\ran}{{\rangle}}

\newtheorem{theo}{Theorem}[section]

\newtheorem{lem}[theo]{Lemma}

\newtheorem{remark}[theo]{Remark}
\newtheorem{definition}[theo]{Definition}
\newtheorem{example}[theo]{Example}
\newtheorem{property}[theo]{Property}

\newenvironment{rem}{\begin{remark} \rm}{\end{remark}}
\newenvironment{pf}{{\it Proof.}}{\hfill$\square$\vspace{3mm}}

\ncm{\R}{\mbox{\bf R}} \ncm{\Z}{\mbox{\bf Z}} \ncm{\T}{\mbox{\bf 
T}} \ncm{\TT}{\T$^{2}$} \ncm{\N}{\mbox{\bf N}} \ncm{\C}{\mbox{\bf 
C}} 

\title{Homogeneity of the pure state space of a  separable \cstar}
\bigskip     
\author{ Akitaka Kishimoto, Narutaka Ozawa, and Sh\^{o}ichir\^{o} Sakai}
\date{}
\begin{document}
\maketitle 

\begin{abstract} We prove that the pure state space is 
homogeneous under the action of the automorphism group (or the 
subgroup of asymptotically inner automorphisms) for all the 
separable simple \cstars. The first result of this kind was shown 
by Powers for the UHF algbras some 30 years ago.  

 \hfill Mathematics Subject 
Classification: 46L40, 46L30 
\end{abstract} 

 \ncm{\U}{{\mathcal{U}}} 
 \ncm{\F}{{\mathcal{F}}}
 \ncm{\G}{{\mathcal{G}}}
 \ncm{\K}{{\mathcal{K}}}
 \ncm{\M}{{\mathcal{M}}}
 \ncm{\E}{{\mathcal{E}}} 
 \ncm{\B}{{\mathcal{B}}}
  \ncm{\vR}{{\mathcal{R}}}   
 \ncm{\om}{{\omega}}
 \ncm{\al}{{\alpha}}
 \ncm{\Hil}{{\mathcal{H}}}   
 \ncm{\Lin}{{\mathcal{L}}}

\section{Introduction}
If $A$ is a \cstar, an automorphism $\al$ of $A$ is {\em 
asymptotically inner} if there is a continuous family 
$(u_t)_{t\in[0,\infty)}$ in the group  $\U(A)$ of unitaries in $A$ 
(or $A+\C1$ if $A$ is non-unital) such that 
$\al=\lim_{t\rightarrow\infty}\Ad\,u_t$; we denote by $\AI(A)$ the 
group of asymptotically inner automorphisms of $A$, which is a 
normal subgroup of the group of approximately inner automorphisms. 
Note that each $\al\in\AI(A)$ leaves each (closed two-sided) ideal 
of $A$ invariant. It is shown, in \cite{Pow,Br,FKK}, for a large 
class of separable \cstars\ that if $\om_1$ and $\om_2$ are pure 
states of $A$ such that the GNS representations associated with 
$\om_1$ and $\om_2$ have the same kernel, then there is an 
$\al\in\AI(A)$ such that $\om_1=\om_2\al$. We shall show in this 
paper that this is the case for {\em all} the separable \cstars; 
formally, denoting by $\pi_{\omega}$ the GNS representation 
associated with a state $\omega$, we state: 

\begin{theo} \label{T}
Let $A$ be a separable \cstar. If $\omega_1$ and $\omega_2$ are 
pure states of $A$ such that $\ker 
\pi_{\omega_1}=\ker\pi_{\omega_2}$, then there is an 
$\alpha\in\AI(A)$ such that $\omega_1\alpha=\omega_2$. 
\end{theo} 

In particular the pure state space of a separable simple  \cstar\ 
$A$ is homogeneous under the action of $\AI(A)$. We need the 
separability for this statement to be true even if we replace 
$\AI(A)$ by the full automorphism group $\Aut(A)$ (see \ref{E}). 
But if we instead assume that $A$ is nuclear, the situation is 
unclear, i.e., we do not know if the pure state space of a 
non-separable simple nuclear \cstar\ is homogeneous under the 
action of $\Aut(A)$ or not. See \cite{Ef} for some problems on 
this.   

We note here that $\AI(A)$ can be considered as a {\em core} of 
$\Aut(A)$ whose inner structure is beyond algebraic grasp; 
$\AI(A)$ is characterized as the subgroup of automorphisms which 
have the same KK class with the identity automorphism for the 
class of purely infinite simple separable \cstars\ classified by 
Kirchberg and Phillips \cite{KP} (see \cite{KK} for a similar 
result for a class of AT algebras). 
                                         
The proof of the above theorem comprises three observations taken 
from \cite{FKK} and \cite{Haa0}. By combining these, the theorem 
will follow immediately.  

The first observation from \cite{FKK} is that the following 
property for a \cstar\ $A$ will imply the above theorem. 

\begin{property} \label{A}
For any finite subset $\F$ of $A$, any pure state $\om$ of $A$ 
with $\pi_{\om}(A)\cap \K(\Hil_{\om})=(0)$, and $\eps>0$, there 
exist a finite subset $\G$ of $A$ and $\delta>0$ satisfying: If 
$\varphi$ is a pure state of $A$ such that $\pi_\varphi$ is 
quasi-equivalent to $\pi_\om$, and 
 $$
 |\varphi(x)-\om(x)|<\delta,\ \ x\in\G,
 $$
then there is a continuous path $(u_t)_{t\in[0,1]}$ in $\U(A)$ 
such that $u_0=1,\ \varphi=\om\Ad\,u_1$, and 
 $$
 \|\Ad\,u_t(x)-x\|<\eps,\ \ x\in\F,\ t\in[0,1].
 $$
\end{property}

In the above statement, $\K(\Hil_{\om})$ is the \cstar\ of compact 
operators on $\Hil_{\om}$, the Hilbert space for $\pi_{\om}$. 
                                                             
Another observation from \cite{FKK} is that the following property 
of $A$, a kind of weak amenability, implies the above property: 

\begin{property}  \label{B}
Let $\F$ be a finite subset of $A$, $\pi$ an irreducible 
representation of $A$ on a Hilbert space $\Hil$, $E$ a 
finite-dimensional projection on $\Hil$, and $\eps>0$. Then there 
exists an $x=(x_1,x_2,\ldots,x_n)\in M_{1n}(A)$ for some $n$ such 
that $\|xx^*\|\leq1$, $\pi(xx^*)E=E$, and $\|\ad\,a\Ad\,x\|<\eps$ 
for all $a\in\F$, where $\Ad\,x$ and $\ad\,a$ denote the linear 
maps on $A$ defined by $b\mapsto xbx^*=\sum x_ibx_i^*$ and 
$b\mapsto  [a,b]$ respectively. 
\end{property}  

Here $M_{mn}(A)$ denotes the $m\times n$ matrices over $A$.
                   
The final observation, from Haagerup \cite{Haa0}, is that this 
property holds for {\em all}  \cstars, which is shown by 
repeating, almost verbatim, the proof of 3.1 of \cite{Haa0} 
employed for verifying the statement that all nuclear \cstars\ are 
amenable.  

Although those observations are mostly immediate from the cited 
references if once properly formulated as above, we shall outline 
the  proofs for the reader's convenience: \ref{A} implies \ref{T} 
in section 2, \ref{B} implies \ref{A} in section 3, and Property 
\ref{B} is universal in section 4. 

The present method is further exploited in connection with 
one-parameter automorphism groups \cite{K01} and for type III 
representations \cite{FKK2}.             
  
\section{Homogeneity} 
 \setcounter{theo}{0}
We denote by $\AI_0(A)$ the set of $\alpha\in\AI(A)$ which has a 
continuous family $(u_t)_{t\in [0,\infty)}$ in $\U(A)$ with 
$u_0=1$ and $\alpha=\lim_{t\rightarrow\infty}\Ad\,u_t$.

\begin{theo}\label{C}
Let $A$ be a separable \cstar\ satisfying Property \ref{A}. If 
$\om_1$ and $\om_2$ are pure states of $A$ such that 
$\ker\pi_{\om_1}=\ker\pi_{\om_2}$, then there is an 
$\al\in\AI_0(A)$ such that $\om_1=\om_2\al$. 
\end{theo}

The following gives a slightly weaker version of Property \ref{A}. 
  
\begin{lem}\label{B1}
Let $A$ be a  \cstar\ with Property \ref{A}. Then for any finite 
subset $\F$ of $A$, any pure state $\om$ of $A$ with 
$\pi_{\om}(A)\cap \K(\Hil_{\om})=(0)$, and $\eps>0$, there exist a 
finite subset $\G$ of $A$ and $\delta>0$ satisfying: If $\varphi$ 
is a pure state of $A$ such that 
$\ker\pi_{\varphi}=\ker\pi_{\om}$, and 
 $$
 |\varphi(x)-\om(x)|<\delta,\ \ x\in\G,
 $$
then for any finite subset $\F'$ of $A$ and $\eps'>0$ there is a 
continuous path $(u_t)_{t\in[0,1]}$ in $\U(A)$ such that $u_0=1$, 
and 
 \BE
 |\varphi(x)-\om\Ad\,u_1(x)|<\eps', \ \ && x\in\F',\\
 \|\Ad\,u_t(x)-x\|<\eps,\ \ &&x\in \F.
 \EE
\end{lem}
\begin{pf}
Given $(\F,\om,\eps)$, choose $(\G,\delta)$ as in Property 
\ref{A}. Let $\varphi$ be a pure state of $A$ such that 
$\ker\pi_{\varphi}=\ker\pi_{\om}$ and 
 $$
 |\varphi(x)-\om(x)|<\delta/2,\ \ x\in\G.
 $$
Let $\F'$ be a finite subset of $A$ and $\eps'>0$ with 
$\eps'<\delta/2$. We can mimic $\varphi$ as a vector state through 
$\pi_{\om}$; by Kadison's transitivity there is a $v\in\U(A)$ such 
that 
 $$
 |\varphi(x)-\om\Ad\,v(x)|<\eps',\ \ x\in\F'\cup\G,
 $$
(see 2.3 of \cite{FKK}). Since $|\om\Ad\,v(x)-\om(x)|<\delta,\ 
x\in\G$, we have, by applying Property \ref{A} to the pair $\om$ 
and $\om\Ad\,v$, a continuous path $(u_t)$ in $\U(A)$ such that 
$u_0=1$, and 
 \BE
 &&\om\Ad\,v=\om\Ad\,u_1,\\
 &&\|\Ad\,u_t(x)-x\|<\eps,\  \ x\in\F.
 \EE
Since $|\varphi(x)-\om\Ad\,u_1(x)|<\eps',\ x\in\F'$, this 
completes the proof.    
\end{pf} 

We shall now turn to the proof of Theorem \ref{C}.
     
Once we have Lemma \ref{B1}, we can prove this in the same way as 
2.5 of \cite{FKK}. We shall only give an outline here. 

Let $\om_1$ and $\om_2$ be pure states of $A$ such that  
$\ker\pi_{\om_1}=\ker\pi_{\om_2}$. 

If $\pi_{\om_1}(A)\cap \K(\Hil_{\om_1})\not=(0)$, then  
$\pi_{\om_1}(A)\supset \K(\Hil_{\om_1})$ and $\pi_{\om_1}$ is 
equivalent to $\pi_{\om_2}$. Then by Kadison's transitivity (see, 
e.g., 1.21.16 of \cite{Sak}), there is a continuous path $(u_t)$ 
in $\U(A)$ such that $u_0=1$ and $\om_1=\om_2\Ad\,u_1$. 
                  
Suppose that  $\pi_{\om_1}(A)\cap \K(\Hil_{\om_1})=(0)$, which 
also implies that $\pi_{\om_2}(A)\cap \K(\Hil_{\om_2})=(0)$. 

Let $(x_n)$ be a dense sequence in $A$. 

Let $\F_1=\{x_1\}$ and $\eps>0$ (or $\eps=1$). Let 
$(\G_1,\delta_1)$ be the $(\G,\delta)$ for $(\F_1,\om_1,\eps/2)$ 
as in Lemma \ref{B1} such that $\G_1\supset\F_1$. For this 
$(\G_1,\delta_1)$ we choose a continuous path $(u_{1t})$ in 
$\U(A)$ such that $u_{1,0}=1$ and 
 $$
 |\om_1(x)-\om_2\Ad\,u_{1,1}(x)|<\delta_1,\ \ x\in\G_1.
 $$
Let $\F_2=\{x_i,\Ad\,u_{1,1}^*(x_i)\ | \ i=1,2\}$ and let 
$(\G_2,\delta_2)$ be the $(\G,\delta)$ for 
$(\F_2,\om_2\Ad\,u_{1,1},2^{-2}\eps)$ as in Lemma \ref{B1} such 
that $\G_2\supset\G_1\cup\F_2$ and $\delta_2<\delta_1$. By 
\ref{B1} there is a continuous path $(u_{2t})$ in $\U(A)$ such 
that $u_{2,0}=1$ and 
 \BE
 \|\Ad\,u_{2t}(x)-x\|<2^{-1}\eps,\ \ &&x\in \F_1,\\
 |\om_2\Ad\,u_{1,1}(x)-\om_1\Ad\,u_{2,1}(x)|<\delta_2,\ \ &&x\in\G_2.
 \EE
Let $\F_3=\{x_i,\Ad\,u_{2,1}^*(x_i)\ | \ i=1,2,3\}$ and let 
$(\G_3,\delta_3)$ be the $(\G,\delta)$ for 
$(\F_3,\om_1\Ad\,u_{2,1},2^{-3}\eps)$ as in \ref{B1} such that 
$\G_3\supset\G_2\cup\F_3$ and $\delta_3<\delta_2$. By \ref{B1} 
there is a continuous path $(u_{3t})$ in $\U(A)$ such that 
$u_{3,0}=1$ and 
 \BE
 \|\Ad\,u_{3t}(x)-x\|<2^{-2}\eps,\ \ &&x\in\F_2,\\    
 |\om_1\Ad\,u_{2,1}(x)-\om_2\Ad(u_{1,1}u_{3,1})(x)|<\delta_3,\ \ 
 &&x\in \G_3.
 \EE
We shall repeat this process. 

Assume that we have constructed $\F_n,\G_n,\delta_n$, and 
$(u_{n,t})$ inductively. In particular if $n$ is even, $\F_n$ is 
given as 
 $$
 \{x_i,\Ad(u_{n-1,1}^*u_{n-3,1}^*\cdots u_{1,1}^*)(x_i)\ |\ i=1,2,\ldots,n\}
 $$
and $(G_n,\delta_n)$ is the $(\G,\delta)$ for 
$(\F_n,\om_2\Ad(u_{1,1}u_{3,1}\cdots u_{n-1,1}),2^{-n}\eps)$ as in 
\ref{B1} such that $\G_n\supset \G_{n-1}\cup\F_n$ and 
$\delta_n<\delta_{n-1}$. And $(u_{n,t})$ is given by \ref{B1} for 
$(\F_{n-1},\om_1\Ad(u_{2,1}\cdots u_{n-2,1}),2^{-n+1}\eps)$ and 
for $\F'=\G_n$ and $\eps'=\delta_n$ and it satisfies 
 \BE
 && \|\Ad\,u_{nt}(x)-x\|<2^{-n+1}\eps, \ \ x\in \F_{n-1},\\
 &&|\om_1\Ad(u_{2,1}u_{4,1}\cdots u_{n,1})(x)
 -\om_2\Ad(u_{1,1}u_{3,1}\cdots u_{n-1,1})(x)|<\delta_n,\ \ x\in\G_n.
 \EE
We define continuous paths $(v_t)$ and $(w_t)$ in $\U(A)$ with 
$t\in[0,\infty)$ by: For $t\in [n,n+1]$ 
 \BE
 &&v_t=u_{1,1}u_{3,1}\cdots u_{2n-1,1}u_{2n+1,t-n},\\
 &&w_t=u_{2,1}u_{4,1}\cdots u_{2n-2,1}u_{2n+2,t-n}.
 \EE
Then, since $\|\Ad\,u_{nt}(x)-x\|<2^{-n+1}\eps,\ x\in \F_{n-1}$, 
we can show that $\Ad\,v_t$ (resp. $\Ad\,w_t$) converges to an 
automorphism $\al$ (resp. $\beta$) as $t\ra\infty$ and that $ 
\om_1\beta=\om_2\al$. Since $\al,\beta\in \AI_0(A)$ and $\AI_0(A)$ 
is a group, this will complete the proof. See the proofs of 2.5 
and 2.8 of \cite{FKK} for details. 

\medskip 

\begin{rem}\label{E}
Let $A$ be a factor of type II$_1$ or type III with separable 
predual $A_*$, which is a unital simple non-separable non-nuclear 
\cstar. Then the pure state space of $A$ is not homogeneous under 
the action of the automorphism group $\Aut(A)$ of $A$. \end{rem} 

This is shown as follows. Since $A$ contains a $C^*$-subalgebra 
isomorphic to $C_b(\N)\equiv C(\beta\N)$ and $\beta\N$ has 
cardinality $2^c$, the pure state space of $A$ has cardinality (at 
least) $2^c$, where $c$ denotes the cardinality of the continuum. 
(We owe this argument to J. Anderson.) On the other hand any 
$\alpha\in\Aut(A)$ corresponds to an isometry on the predual 
$A_*$, a separable Banach space.  Thus, since the set of bounded 
operators on a separable Banach space has cardinality $c$, 
$\Aut(A)$ has cardinality (at most) $c$. Hence  the pure state 
space of $A$ cannot be homogeneous under the action of  $\Aut(A)$.

\section{\ref{B} implies \ref{A}}   
 \setcounter{theo}{0}                                                               

\begin{theo}
Any \cstar\ with Property \ref{B} has  Property \ref{A}. 
\end{theo} 
\begin{pf}
Let $\F$ be a finite subset of $A$, $\om$ a pure state of $A$ with 
$\pi_{\om}(A)\cap \K(\Hil_{\om})=(0)$, and $\eps>0$.  For 
$\pi=\pi_{\om}$ and the projection $E$ onto the subspace 
$\C\Omega_{\om}$, we choose an $x\in M_{1n}(A)$ for some $n$ as in 
Property \ref{B}, i.e., $\|x\|\leq1$, 
$\pi(xx^*)\Omega_{\om}=\Omega_{\om}$ with 
$\Omega=\Omega_{\omega}$, and $\|\ad\,a\,\Ad\,x\|<\eps$ for all 
$a\in \F$. 

Let  
 $$
 \G=\{x_ix_j^*\ |\ i,j=1,2,\ldots,n\},
 $$                                   
which will be  the subset $\G$ required in Property \ref{A}. We 
will choose $\delta>0$ sufficiently small later. Suppose that we 
are given a unit vector $\eta\in\Hil_{\om}$ satisfying 
 $$
 |\lan\pi(x_i^*)\eta,\pi(x_j^*)\eta\ran-\lan\pi(x_i^*)\Omega,\pi(x_j^*)\Omega\ran|<\delta
 $$
for  any $i,j=1,2,\ldots,n$, where $\Omega=\Omega_{\om}$. Note 
that 
 $$
 \sum_{j=1}^n\|\pi(x_j^*)\Omega\|^2=\lan\pi(xx^*)\Omega,\Omega\ran=1,
 $$
which implies, in particular, that 
$|\lan\pi(xx^*)\eta,\eta\ran-1|<n\delta$.  Thus the two finite 
sets of vectors $S_{\Omega}=\{\pi(x_i^*)\Omega\ |\ i=1,\ldots,n 
\}$ and $S_{\eta}=\{\pi(x_i^*)\eta\ |\ i=1,\ldots,n \}$ have 
similar geometric properties in $\Hil_{\omega}$ if $\delta$ is 
sufficiently small.  Hence we are in a situation where we can 
apply 3.3 of \cite{FKK}. 

Let us describe how we proceed from here in a simplified case. 
Suppose that the linear span $\Lin_{\Omega}$ of $S_{\Omega}$ is 
orthogonal to the linear span $\Lin_{\eta}$ of $S_{\eta}$ and that 
the map $\pi(x_i^*)\Omega\mapsto \pi(x_i^*)\eta$ and 
$\pi(x_i^*)\eta\mapsto  \pi(x_i^*)\Omega$ extends to a unitary $U$ 
on $\Lin_{\Omega}+\Lin_{\eta}$; in particular we have assumed that 
$\lan\pi(x_i^*)\eta,\pi(x_j^*)\eta\ran 
=\lan\pi(x_i^*)\Omega,\pi(x_j^*)\Omega\ran$ for all $i,j$. Since 
$U$ is a self-adjoint unitary, $F\equiv(1-U)/2$ is a projection 
and satisfies that $e^{i\pi F}=U$ on the finite-dimensional 
subspace $\Lin_{\Omega}+\Lin_{\eta}$. By Kadison's transitivity we 
choose an $h\in A$ such that $0\leq h\leq 1$ and 
$\pi(h)|\Lin_{\Omega}+\Lin_{\eta}=F$. We set $ \ol{h}=\Ad\,x(h)$, 
which entails that $\|[a,\ol{h}]\|<\eps,\ a\in\F$. Then we have 
that  
 \BE 
 \pi(\ol{h})(\Omega-\eta)&=&
 \pi(xhx^*)(\Omega-\eta)\\&=&\sum\pi(x_i)F\pi(x_i^*)(\Omega-\eta),\\ 
 &=&\sum\pi(x_i)\pi(x_i^*)(\Omega-\eta)\\
 &=&\Omega-\eta
 \EE
and that $\pi(\ol{h})(\Omega+\eta)=0$. Hence it follows that
 $$
 \pi(e^{i\pi \ol{h}})\Omega=\pi(e^{i\pi \ol{h}})(\Omega-\eta)/2
 +\pi(e^{i\pi 
 \ol{h}})(\Omega+\eta)/2=-(\Omega-\eta)/2+(\Omega+\eta)/2=\eta.
 $$

Thus the path $(e^{it\pi \ol{h}})_{t\in[0,1]}$  is what is 
desired. 
 
Whenever $\calL_{\Omega}$ is orthogonal to $\calL_{\eta}$, this 
argument can be made rigorous if $\delta>0$ is sufficiently small. 
See \cite{FKK} for details. 

If $\Lin_{\eta}$ is not orthogonal to $\Lin_{\Omega}$, we still 
find a unit vector $\zeta\in\Hil_{\om}$ such that 
  $$
 |\lan\pi(x_i^*)\zeta,\pi(x_j^*)\zeta\ran-\lan\pi(x_i^*)\Omega,\pi(x_j^*)\Omega\ran|<\delta
 $$
and such that $\Lin_{\zeta}$ is orthogonal to both $\Lin_{\Omega}$ 
and $\Lin_{\eta}$. Here we use the assumption that 
$\pi_{\om}(A)\cap\K(\Hil_{\om})=(0)$. Then we combine the path of 
unitaries sending $\eta$ to $\zeta$ and then the path sending 
$\zeta$ to $\Omega$ to obtain the desired path. 
\end{pf}

\section{Property \ref{B} is universal} 
 \setcounter{theo}{0}    
\ncm{\Bil}{{\rm Bil}} \newcommand{\ip}[1]{\langle#1\rangle} 
\newcommand{\e}{\varepsilon}

Let $\Bil(A)$ denote the bounded bilinear forms on a \cstar\ $A$. 
We have a canonical isometric identification of $\Bil(A)$ with 
$(A\widehat{\otimes}A)^*$, which is given by 
$$\ip{V,a\otimes b}=V(a,b).$$                             
Here $A\widehat{\otimes}A$ is the completion of the algebraic 
tensor product $A\otimes A$ equipped with the projective tensor 
norm: 
$$\| S\|_{\wedge}=\inf\{\sum_{i=1}^n\| x_i\|\ \| y_i\|\},
$$ 
where the infimum is taken all over the possible representations 
$S=\sum_{i=1}^n x_i\otimes y_i$. For $a\in A$ the bounded linear 
maps $L_a$ and $R_a$ on $A\widehat{\otimes}A$ are defined by 
$$L_a(x\otimes y)=ax\otimes y\quad
\mbox{ and }\quad R_a(x\otimes y)=x\otimes ya$$ and the bounded 
linear map $p\colon A\widehat{\otimes}A\to A$ is defined by 
$$p(x\otimes y)=xy.$$

If $\M$ is a von Neumann algebra, $\Bil_{\sigma}(\M)$ denotes the 
subspace of $\Bil(\M)$ consisting of separately $\sigma$-weakly 
continuous forms on $\M$. For $a\in\M$, the dual maps $(L_a)^*$ 
and $(R_a)^*$ leave $\Bil_{\sigma}(\M)$ invariant. We define a 
contraction $\varphi: \Bil(\M)\ra \ell^{\infty}(\M_1)$ by 
$\varphi(V)(a)=V(a^*,a)$, where $\M_1$ is the unit ball of $\M$.

 We rely on the following result \cite{Haa0}: 

\begin{theo} (Haagerup)
Let $\M$ be an injective von Neumann algebra. Then there exists a 
mean $m$ on the (discrete) semigroup $I(\M)$ of isometries in $\M$ 
which is invariant in the sense that 
 $$
 m(\varphi(L_a^*V)|I(\M))=m(\varphi(R_a^*V)|I(\M))
 $$
for all $V\in \Bil_{\sigma}(\M)$ and all $a\in \M$. \end{theo} 

By using the above result and the proof of 3.1 of \cite{Haa0} we 
prove:

\begin{lem}
Let $\pi\colon A\to \B(\Hil)$ be a non-degenerate representation 
of a $C^*$-algebra $A$. If $\pi(A)''$ is injective, then there 
exists a net $\{T_\lambda\}_\lambda$ in $A\otimes A$ such that 
\begin{enumerate}
\item the net $\{ T_\lambda\}$ is in the convex hull of 
$\{ x\otimes x^*\, |\ x\in A,\ \| x\|\le 1\}$, 
\item $\lim_{\lambda}\| L_aT_\lambda-R_aT_\lambda\|_{\wedge}=0$ 
for any $a\in A$, 
\item $\pi(p(T_\lambda))\to 1$ $\sigma$-weakly in $\B(\Hil)$.
\end{enumerate} 
\end{lem}
\begin{pf}                       
What is shown as Theorem 3.1 in \cite{Haa0} is the above statement 
(or more precisely the statement on $\omega$ below) for a nuclear 
\cstar\ $A$ and its universal representation $\pi$. But the proof 
there depends only on the fact that $\M=\pi(A)''$ is injective. We 
shall just give an outline of the proof here. 

Let $e$ denote the central projection in $A^{**}$ corresponding to 
$\pi$; we shall identify $\M$ with $A^{**}e$. 

By using the fact that $V\in\Bil(A)$ uniquely extends to 
$\tilde{V}\in\Bil_{\sigma}(A^{**})$ \cite{JKR}, we define an 
$\omega\in (A\widehat{\otimes}A)^{**}\cong \Bil(A)^*$ by 
 $$
 \omega(V)=m(\varphi(\tilde{V})|I(\M)),
 $$                                    
where $m$ is an {\em invariant} mean on $I(\M)$ as in the above 
theorem. We then assert that 
\begin{enumerate} 
\item $\omega$ is in the weak$^*$-closed convex hull of 
$\{ x\otimes x^*\, |\ x\in A,\ \| x\|\le1\}$, 
\item $L_a^{**}\omega=R_a^{**}\omega$ for any $a\in A$, 
\item $p^{**}(\omega)=e$ in $A^{**}$, 
\end{enumerate}                      
Property 1 follows by the Hahn-Banach separation argument using 
the crucial fact that $\tilde{V}$ is jointly $\sigma$-strong$^*$ 
continuous \cite{Haa1}. Property 2 reflects the invariance of $m$ 
in the above theorem: 
$(L_a^{**}\omega)(V)=\omega(L_a^*V)=m(\varphi(L_a^*\tilde{V})|I(\M)) 
=m(\varphi(L_{ae}^*(\tilde{V}|\M)|I(\M)))=m(\varphi(R_{ae}^*(\tilde{V}|\M)|I(\M)))$, 
which is equal to $(R_a^{**}\omega)(V)$, for all $V\in\Bil(A)$ and 
$a\in A$, where $\tilde{V}|\M\in \Bil_{\sigma}(\M)$ is the 
restriction of $\tilde{V}$. Since $(p^*f)\tilde{}=p^*f$ and 
$\varphi(p^*f)(a)=f(a^*a)$ for $f\in A^*$, Property 3 follows 
from: $p^{**}(\omega)(f)=\omega(p^*f)=f(e)$. 

Now, we may find a net $\{ T_\lambda\}$ in the convex hull of $\{ 
x\otimes x^*\, |\ x\in A,\ \| x\|\le1\}$ such that $T_\lambda$ 
weak$^*$-converges to $\omega$ in $(A\widehat{\otimes}A)^{**}$. It 
follows that $p(T_\lambda)$ weak$^*$-converges to $e$ in $A^{**}$. 
Since for any $a\in A$, $L_aT_\lambda-R_aT_\lambda$ converges 
weakly to $0$ in $A\widehat{\otimes}A$, we may assume 
$$\lim_\lambda\| L_aT_\lambda-R_aT_\lambda\|_{\wedge}=0$$ 
by convexity. 
\end{pf}
         
By applying the above lemma to an irreducible representation $\pi$ 
on $\Hil$, a finite-dimensional projection $E$ on $\Hil$, and 
$\eps>0$,  we obtain a sequence $(x_1,x_2,\ldots,x_n)$ in $A$ such 
that $\sum_{i=1}^n\|x_i\|^2\leq 1$, and 
 \BE
 &&\|\sum_iax_i\otimes x_i^*-\sum_ix_i\otimes x_i^*a\|_\wedge<\eps,\ \ a\in\F,\\
 &&\|\pi(\sum_ix_ix_i^*)E-E\|<\eps.
 \EE
By using Kadison's transitivity, we find a $b\in A\, ({\rm or}\ 
A+\C1)$ such that $b\approx 1$, $\|yy^*\|\leq 1$, and 
$\pi(yy^*)E=E$, where $y=(bx_1,bx_2,\ldots,bx_n)\in M_{1n}(A)$.  
Since there is a contraction $\psi$ of $A\widehat{\otimes}A$ into 
$\B(A)$, which is defined by $\psi(a\otimes b)(x)=axb$, we obtain:

\begin{theo}
Any \cstar\ has Property \ref{B}. \end{theo}

\medskip
\small
\begin{flushright}
 Department of Mathematics, Hokkaido University, Sapporo, Japan 060-0810 \\ 
 Department of Mathematical Sciences, University of Tokyo, Tokyo, Japan 153-8914\\

 5-1-6-205, Odawara, Aoba-ku,  Sendai, Japan 980-0003 
\end{flushright}


\small
\begin{thebibliography}{99}
\bibitem{Br} O. Bratteli, Inductive limits of finite-dimensional 
\cstars, Trans. Amer. Math. Soc. 171 (1972), 195--234. 
\bibitem{Ef} E.G. Effros, On the structure theory of \cstars: 
some old and new problems, in: Proceedings of symposia in pure 
mathematics 38 (1982) part 1, edited by R.V. Kadison, pages 
19--34. 
\bibitem{FKK} H. Futamura, N. Kataoka, and A. Kishimoto, Homogeneity 
of the pure state space for separable \cstars, to appear in 
Internat. J. Math. 
\bibitem{FKK2}  H. Futamura, N. Kataoka, and A. Kishimoto, Type 
III representations and automorphisms of some separable nuclear 
\cstars, in preparation. 

 
\bibitem{Haa0} U. Haagerup, All nuclear \cstars\ are amenable, Invent. Math. 74 (1983),
305--319.    
\bibitem{Haa1} U. Haagerup, The Grothendieck inequality for 
bilinear forms on \cstars, Adv. in Math. 56 (1985), 93--116.
 
\bibitem{KP} E. Kirchberg and N.C. Phillips, Embedding of exact \cstars\ and continuous
 fields in the Cuntz algebra ${\cal O}_2$, J. reine angew. Math.  525 (2000), 
 55--94.
\bibitem{K01} A. Kishimoto, Approximately inner flows on separable 
\cstars, preprint.
\bibitem{KK}  A. Kishimoto and A. Kumjian, The Ext class of 
approximately inner automorphisms, II, to appear in J. Operator 
Theory.                                                    
\bibitem{JKR} B.E. Johnson, R.V. Kadison, and J.R. Ringrose, 
Cohomology of operator algebras III. Reduction to normal 
cohomology, Bull. Soc. Math. France 100 (1972), 73--96. 
 \bibitem{Pow} R.T. Powers, Representations of uniformly hyperfinite 
algebras and their associated von Neumann rings, Ann. of Math. 86 
(1967), 138--171. 
\bibitem{Sak} S. Sakai, {\em $C^*$-algebras and $W^*$-algebras}, 
Classics in Math., Springer, 1998. 
\end{thebibliography}
\end{document}